\documentclass[12pt]{article}
\usepackage{amssymb, cite}

\newcommand{\be}{\begin{equation}}
\newcommand{\ee}{\end{equation}}
\newcommand{\bqn}{\begin{eqnarray}}

\newcommand{\eqn}{\end{eqnarray}}

\newcommand{\bd}{\begin{description}}
\newcommand{\ed}{\end{description}}

\newtheorem{stat}{}[section]

\def\bs{\begin{stat}}
\def\es{\end{stat}}
\def\ben{\begin{enumerate}}
\def\een{\end{enumerate}}

\def\bp{\noindent{\bf Proof}  \ }
\newcommand{\ep}{\hfill $\square$}

\makeatletter
\@addtoreset{equation}{section}

\makeatother

\begin{document}

\begin{center}
{\bf COUNTEREXAMPLES 
\\[0.7ex]
TO THE  CUBIC GRAPH DOMINATION CONJECTURE}
\\[4ex]
{\large {\bf Alexander Kelmans}}
\\[2ex]
{\bf University of Puerto Rico, San Juan, Puerto Rico}
\\[0.5ex]
{\bf Rutgers University, New Brunswick, New Jersey}
\\[2ex]
\end{center}

\begin{abstract}
Let $v(G)$ and $\gamma (G)$ denote the number of vertices and the domination number of a graph $G$, respectively, and
let $\rho (G) = \gamma (G)/v(G)$.
In 1996 B. Reed 
conjectured that if $G$ is  a cubic graph, then $\gamma (G) \le \lceil v(G)/3 \rceil $. In 2005 A. Kostochka and B. Stodolsky disproved this conjecture for cubic graphs of connectivity one and
 maintained that the conjecture may still be true for cubic 
 2-connected graphs.
Their  minimum counterexample $C$ has 4 bridges, 
$v(C) = 60$, and $\gamma (C) = 21$.
In this paper we disprove Reed's conjecture for cubic 
2-connected graphs by providing a sequence 
$(R_k:  k \ge 3)$ of  cubic graphs of connectivity two
with 
$\rho (R_k) = \frac{1}{3} + \frac{1}{60}$,
where $v(R_{k+1}) > v(R_k) > v(R_3) = 60$ for $k \ge 4$, and so 
$\gamma (R_3) = 21$ and 
$\gamma (R_k) -  \lceil v(R_k)/3 \rceil \to \infty $ with $k \to \infty $.
%
%
We also provide a sequence of $(L_r: r \ge 1)$ 
of  cubic graphs of connectivity one
with $\rho (L_r) > \frac{1}{3} + \frac{1}{60}$.
The minimum counterexample $L = L_1$ in this sequence
is `better' than $C$ in the sense that
$L$ has 2 bridges while $C$ has 4 bridges, 
$v(L) =  54 < 60 = v(C)$, and
$ \rho (L) = \frac{1}{3} + \frac{1}{54} > 
\frac{1}{3} + \frac{1}{60} = \rho (C)$.
We also give a construction providing for every $r \in \{0,1,2\}$ infinitely many cubic cyclically 4-connected
Hamiltonian graphs $G_r$ such that $v(G_r) = r \bmod 3$,
$r \in \{0,2\} \Rightarrow \gamma (G_r) = \lceil v(G_r)/3 \rceil $,
and 
$r = 1 \Rightarrow \gamma (G_r) = \lfloor v(G_r)/3 \rfloor $.
At last we suggest a stronger conjecture on domination in cubic 3-connected graphs.
\\[1ex]
\end{abstract}


\section{Introduction}

\indent

We consider simple undirected graphs. 
All notions on graphs that are  not defined here can be found in \cite{D}.

Let $G$ be a graph, $V(G)$ and $E(G)$ the sets of vertices and  edges of  $G$, respectively,
$v(G) = |V(G)|$ and $e(G) = |E(G)|$.
Let $N(v,G)$ denote the set of vertices in $G$ adjacent to 
a vertex $v$. Let $\kappa (G)$ denote the vertex connectivity
of $G$.
A vertex subset $X$ of $G$ is called {\em dominating} if every vertex in $G - X$ is adjacent to a vertex in $X$.
Let $\gamma (G)$ denote   the size of a minimum dominating set in $G$;  $\gamma (G)$ is  called the {\em dominating number}  of $G$. 
We call $\rho (G) = \gamma (G)/ v(G)$ the {\em dominating ratio} of $G$.
A graph $G$ is called {\em cubic} if every vertex of $G$ has degree three.

Quite a few papers (e.g. \cite{Ar,BC,CH,KPS,MS,KS,O,R}), 
a survey paper 
\cite{CH}, and 
a book \cite{HHS} are devoted to various problems related to
the  domination number  and 
its relations with some other parameters of graphs.


In 1996 \cite{R}, B. Reed proved that if 
the minimum vertex degree in $G$ is  at least three, then 
$\gamma (G) \le 3v(G)/8$ 
and conjectured that if in addition $G$ is cubic,  
then  $\gamma (G) \le \lceil v(G)/3\rceil $.
In 2005 \cite{KS} A. Kostochka and B. Stodolsky disproved  Reed's conjecture for cubic graphs of connectivity one by presenting a sequence of cubic graphs $G$ of connectivity one with
$\rho (G) >  \frac{1}{3} + \frac{1}{69} $ and maintained that Reed's conjecture may still be true for cubic 2--connected graphs.  
Let $C$ and $H$ be the minimum counterexample and  another counterexample 
in \cite{KS}, respectively.
Then $C$ has four bridges,
$v(C) = 60$, and $\rho (C) = \frac{7}{20} = 
\frac{1}{3} + \frac{1}{60} > \rho (H) >  \frac{1}{3} + \frac{1}{69}$.

In this paper we disprove Reed's conjecture for cubic 
2-connected graphs by giving
several constructions 
(see {\bf \ref{R}}, {\bf \ref{gG(P)}}, and {\bf \ref{eRplPi}})  that provide infinitely many  counterexamples  of connectivity two. 
One of our constructions (see {\bf \ref{R}}) provides a sequence $(R_k: k \ge 3)$ of  cubic graphs of connectivity two
with 
$\rho (R_k) = \frac{1}{3} + \frac{1}{60}$,
where $v(R_{k+1}) > v(R_k) > v(R_3) = 60$ for $k \ge 4$, and so 
$\gamma (R_3) = 21$ and 
$\gamma (R_k) -  v(R_k)/3 \to \infty $ with $k \to \infty $.
Thus the violation $\gamma (G) - \lceil v(G)/3 \rceil $
of the inequality   in the Reed's conjecture may be arbitrarily large.
Graph $R_3$ is the minimum 2-connected counterexample we have found. 

We also present (see {\bf \ref{Lr}}) a sequence $(L_r: r \ge 1)$ of  `better' counterexamples of connectivity one than those in \cite{KS}. Namely, 
$L_1$ has two bridges, $v(L_1) = 54$, $v(L_r) < 
\\[0.5ex]
v(L_{r+1})$, and  
$\rho (L_r) = \frac{7}{20} + \frac{1}{200r +340} \to \frac{7}{20}$ with $r \to \infty$,  and so 
$\rho (C) = \frac{1}{3} + \frac{1}{60} < 
\rho (R_k)  < 
\\[0.7ex]
\rho (L_{r+1}) < \rho (L_1) = 
\frac{1}{3} + \frac{1}{54}$. 
Therefore every counterexample in this construction has larger domination ratio than  every counterexample in \cite{KS}. Moreover, $L_1$ has less vertices, 
larger domination ratio, and  less bridges than $C$.


We give constructions 
(see {\bf \ref{vRplB}} and {\bf \ref{Mk}}) that  for every $r \in \{0,1,2\}$ provide infinitely many cubic 3-connected and
cyclically 4-connected graphs $G_r$ such that 
$v(G_r) = r \bmod 3$,
$r \in \{0,2\} \Rightarrow \gamma (G_r) = \lceil v(G_r)/3 \rceil $,
and 
$r = 1 \Rightarrow \gamma (G_r) = \lfloor v(G_r)/3 \rfloor $.


At last we suggest a stronger conjecture (see {\bf \ref{conjecture}}) on domination in cubic 3-connected graphs. 
\\

The results of this paper were discussed in the Department of Mathematics, UPR, in February 2006.

\section{Constructions of counterexamples}

\indent

We start with the following easy observation.

\bs 
\label{HinG} 
Let $G$ be a graph, $H$ an induced subgraph of $G$, and
$X$ the set of vertices in $H$ adjacent to some vertices in 
$G - V(H)$. 
Suppose that $\gamma (H - V) = \gamma (H)$ for every 
$V \subseteq X$. 
If $D$ is a 
dominating set of $G$, then 
$|D \cap V(H)| \ge \gamma (H)$.
\es

Let $H$ be a graph, $\{h_1,h_2\} \subseteq V(H)$, and
$\dot{H} = (H, \{h_1,h_2\})$.
Let $G$ and $H$ be  disjoint graphs and  
$e = v_1v_2 \in E(G)$.
If $G'$ is obtained from $G - e$ and $H$ by identifying 
$h_1$ with $v_1$ and $h_2$ with $v_2$, then we say that 
{\em $G'$  is obtained from $G$ by replacing edge $e$ by 
$\dot{H}$}. 

Let $U$ be a graph, $\{u_1,u_2,u_3\} \subseteq V(U)$, and
$\dot{U} = (U, \{u_1,u_2,u_3\})$.
Let $G$ and $U$ be  disjoint graphs, $v \in V(G)$, and 
$N(v,G) = \{v_1,v_2,v_3\}$.
If $G'$ is obtained from $G - v$ and $U$ by adding three new edges $u_iv_i$, $i \in   \{1,2,3\}$,  then we say that  
{\em $G'$  is obtained from $G$ by replacing vertex $v$ 
by  $\dot{U}$}. 


Let  
$(X, \{x_1, x_2\}$ and $(Y, \{y_1, y_2\})$ be two disjoint copies of $(H, \{h_1,h_2\})$.  Let
$F'$ ($F''$) be obtained from 
$X \cup Y \cup \{x_1y_1, x_2y_2\}$ 
by subdividing edge $x_1y_1$ with a new vertex $z_1$ 
(respectively, by subdividing each edge $x_iy_i$ with a new vertex $z_i$, $i \in \{1,2\}$).

%
%
Let $F_2$ be the graph obtained from $F'' \cup z_1z_2$ by subdividing two  edges $x_1z_1$ and $y_1z_1$ with new vertices $x$ and $y$, respectively. 
Let $F_3$ be the graph obtained from $F_2$ by subdividing
edge $z_1z_2$ with a new vertex $z$. 
Let ${\cal T}_1(\dot{H}) = (F', z_1)$, 
${\cal T}_2(\dot{H}) = (F'', \{z_1, z_2\}))$, 
${\cal F}_2(\dot{H}) = (F_2, \{x, y\}))$, 
and
${\cal F}_3(\dot{H}) = (F_3, \{x, y, z\}))$. 

Let $e_1$, $e_2$, and $e_3$ be three edges in $K_{3,3}$ incident to the same vertex. Let $A$ be the graph obtained from $K_{3,3}$ by subdividing  $e_i$ with 
a new vertex $a_i$  for every $i \in \{1,2\}$.
Similarly, let $B$ be the graph obtained from $K_{3,3}$ 
by subdividing  $e_i$ with  a new vertex $b_i$ for  every 
$i \in \{1,2, 3\}$. 

Let $\dot{A} = (A, \{a_1, a_2\})$ and 
$\dot{B} = (B, \{b_1, b_2, b_3\})$.
Let ${\cal T}_1(\dot{A}) = \dot{S} =(S, s)$ and
${\cal T}_2(\dot{A}) = \dot{T} = (T, \{t_1,t_2\})$
with $\dot{H}: = \dot{A}$).  
Let ${\cal F}_2(\dot{A}) = \dot{P} = (P, \{p_1,p_2\})$ and
${\cal F}_3(\dot{A}) = \dot{Q} = (Q, \{q_1,q_2, q_3\})$
with $\dot{H}: = \dot{A}$).
%
%
\\[2ex]
\indent
It is easy to see the following.

\bs {\em \cite{KS}}
\label{A}
$v(A) = 8$, 
$\gamma (A) = \gamma (A - a_i) = 3$ for every 
$i \in \{1,2\}$, and
$\gamma (A - \{a_1, a_2\}) = 2$.
\es

It is also easy to see the following.

\bs 
\label{B} 
$v(B) = 9$ and  
$\gamma (B - V) =  3$ for every 
$V \subseteq \{b_1, b_2,b_3\}$.
\es

From {\bf \ref{A}} we have:

\bs 
\label{PQS} 
Obviously $v(S) = 17$, $v(T) = 18$, $v(P) = 20$, and 
$v(Q) = 21$. Moreover,
\\[0.5ex]
$(a1)$
$\gamma (S) = \gamma (S - s) = \gamma (T) = 
\gamma (T - t_1) = \gamma (T - t_2) = 
\gamma (T - \{t_1,t_2\}) = 6$,
\\[0.5ex]
$(a2)$
$\gamma (P) = \gamma (P - p_1) = \gamma (P - p_2) = 
\gamma (P - \{p_1,p_2\}) = 7$, and
\\[0.5ex]
$(a3)$ $\gamma (Q - V) = 7$ for every $V \subseteq \{q_1,q_2,q_3\}$.
\es

Let $R_k$ be a graph obtained from a $2k$-vertex cycle 
$(v_0, \ldots , v_{2k-1}, v_{2k})$ with $v_{2k} = v_0$  by replacing each edge $v_{2i}v_{2i+1}$ by a copy  
$(P_i, \{p_1^i, p_2^i\}$ of  $(P,\{p_1,p_2\})$.
%
%

\bs
\label{R} Let $k \ge 3$.
Then $R_k$ is a  cubic graph, $\kappa (R_k) = 2$, 
$v(R_k) = 20k$, and 
$\gamma (R_k) = 7k$, 
\\[0.5ex]
and so 
$\rho (R_k) = \frac{7}{20} = \frac{1}{3} + \frac{1}{60}$ and
$\gamma (R_{k}) - v(R_{k})/3 = k/3 \to \infty $ with 
$k \to \infty $.
\es 

\bp $~$Since $v(P) = 20$, clearly $v(R_k) = 20k$. 
By  {\bf \ref{HinG}} and {\bf \ref{PQS}} $(a2)$, 
$\gamma (R_k) = 7k$.
\ep
\\



Let $T_r$ be obtained from a $2r$-vertex path 
$(v_1, \ldots , v_{2r})$ by replacing each edge $v_{2i-1}v_{2i}$ by a copy  $(P_i,  \{p_1^i, p_2^i\}$ of  $(P, \{p_1,p_2\})$.
Let $L_r = T_r \cup S_1 \cup S_2 \cup \{s_1v_1, s_2v_{2r}\}$, where $(S_1, s_1)$ and  $(S_2, s_2)$ are two copies of $(S,s)$ and  $T_r$, $S_1$, and $S_2$ are disjoint.
\\[2ex]
\indent
From {\bf \ref{HinG}} and {\bf \ref{PQS}} $(a1)$,$(a2)$ 
we have:
\bs
\label{Lr} Let $r \ge 1$.
Then $L_r$ is a cubic graph, $L_r$ has exactly $r+1$
bridges $($and so 
\\[0.5ex]
$\kappa (L_r) = 1$$)$
$v(L_r ) = 20r +34$, and $
\gamma (L_r) = 7r +12$, and so 
$\rho (L_r) = \frac{7}{20} + \frac{1}{200r +340} \to \frac{7}{20}$
\\[0.5ex]
 with $r \to \infty$ and
$\rho (C) = \frac{1}{3} + \frac{1}{60}< \rho (L_{r+1})  < 
\rho (L_r) \le \rho (L_1) = \frac{1}{3} + \frac{1}{54}$.
\es 

Let $P'$ be the graph obtained from $P$ by adding two new vertices $p'_1$, $p'_2$ and two new edges $p_1p'_1$, $p_2p'_2$ and let $\dot{P}' = (P', \{p'_1,p'_2\})$. 
Let 
$G(P)$ be a graph obtained from a  graph $G$ by replacing each edge $e$ by a copy  $\dot{P}'_e$ of $\dot{P}'$.
\\[2ex]
\indent
From {\bf \ref{HinG}} and {\bf \ref{PQS}} $(a2)$ we have:
\bs
\label{G(P)} 
Let $G$ be a  graph.
If $\kappa (G) = 1$, then also $\kappa (G(P)) = 1$.
If $G$ is 2--connected, then $\kappa (G(P)) = 2$.
Also  $v(G(P)) = v(G) + 20e(G)$ and $ \gamma (G(P)) = 7e(G)$.
\es 

From {\bf \ref{G(P)}} we have:
\bs
\label{gG(P)} 
Let $G$ be a connected cubic graph with $2k$ vertices and possible parallel edges.
Then  $v(G(P)) = 62k$,
$\gamma (G(P)) = 21k$, and so
$\rho (G(P)) = \frac{1}{3} + \frac{1}{186}$.
If $\kappa (G) = 1$, then also $\kappa (G(P)) = 1$.
If $G$ is 2--connected, then $\kappa (G(P)) = 2$.
\es 

Given a cubic graph $G$, let $G(P,B)$ be a graph obtained from $G$ by replacing each  vertex $v$ of $G$ by a copy $\dot{B}_v$ of $\dot{B}$ and  each edge $e$ of $G$ 
by a copy  $\dot{P}'_e$ of $\dot{P}'$.
\\[2ex]
\indent
From {\bf \ref{HinG}},  {\bf \ref{B}}, and 
{\bf \ref{PQS}} $(a2)$ we have:

\bs
\label{gG(P,B)} 
Let $G$ be a cubic graph with possible parallel edges and 
with $2k$ vertices. Let $G' = G(P,B)$.
Then   $v(G')  = 78k$, 
$\gamma (G')  = 27k$, and so
$\rho (G')  = \frac{1}{3} + \frac{1}{78}$.
If $\kappa (G) = 1$, then also $\kappa (G') = 1$.
If $G$ is 2--connected, then $\kappa (G') = 2$.
\es



Let us define  $\dot{P}^i$ recursively.
Let $\dot{P}^1 = \dot{P}$  
and
$\dot{P}^{i+1} = {\cal F}_2 (\dot{P}^i)$.
Let ${\cal P} = \{\dot{P}^i: i \ge 1\}$. 

\bs
\label{Pi} 
Let 
$\dot{P}^i = (P^i, \{p_1, p_2\})$, $i \ge 1$.
Then 
\\[0.5ex]
$(a)$ 
$\gamma (P^{i+1}) = 2\gamma (P^i) +1$ and
$\gamma (P^i) = \gamma (P^i - p_1) =  
\gamma (P^i - p_2) = \gamma (P^i - \{p_1, p_2\})$ and 
\\[0.5ex]
$(b)$
$v(P^i) =  2^{i+2} 3  - 4$ and 
$\gamma (P^i) = 2^{i+2}  -1$, and so 
$\rho (P^i) = \frac{1}{3} + \frac{1}{12(2^i3 - 1)}$.
\es

\bp (uses {\bf \ref{PQS}}). Claim $(a)$ can be easily proved 
by induction using 
{\bf \ref{PQS}}.
We prove $(b)$. Obviously 
$v(P^1) = 20$ and by {\bf \ref{PQS}}, $\gamma (P^1) = 7$.
By the definition of $\dot{P}^i$, 
$v(P^{i+1}) = 2v(P^i) +4$.
Now $(b)$ follows from the above recursions for $v(P^{i+1})$ and $\gamma (P^{i+1})$. 
\ep 
\\

Let  $\dot{Q}^i = {\cal F}_3(P^i)$
and
${\cal Q} = \{\dot{Q}^i: i \ge 1\}$. 
From {\bf \ref{PQS}} $(a3)$ and {\bf \ref{Pi}} we have: 
\bs
\label{Ci} 
Let 
$\dot{Q}^i = (Q^i, \{q_1, q_2, q_3\})$.
Then 
\\[0.5ex]
$(a)$
$\gamma (Q^i) = \gamma (Q^i - V)$ for every 
$V  \subseteq  \{q_1, q_2, q_3\}$
and 
\\[0.5ex]
$(b)$
$v(Q^i) =  v(P^i) +1 = 3(2^{i+2}  - 1)$ and 
$\gamma (Q^i) = 2^{i+2}  -1$, and so 
$v(Q^i) = 3 \gamma (Q^i)$.
\es

From {\bf \ref{Pi}} and {\bf \ref{Ci}} we have:
\bs
\label{eRplPi} 
Let $G$ be  either $R_k$ or $L_r$ or
$H(P)$ or $H(P,B)$ for some connected cubic graph $H$. 
Let $G'$ be obtained from $G$ by replacing some copies of $\dot{P}$ 
and/or $\dot{Q}$ in $G$ by copies of some members of 
${\cal P}$ and some copies of $\dot{B}$ by some copies of members of ${\cal Q}$. Then $G'$ is a cubic graph, 
$\gamma (G') > \lceil v(G')/3 \rceil $, and
if $G$ is 2-connected, then $G'$ is also 2-connected.
\es

\section{Cubic 3-connected graphs $G$ with 
$\gamma (G) = \lceil v(G)/3 \rceil $}

\indent


Let $G$ be a cubic graph and 
$G[\dot{B}]$ be a graph obtained from $G$ by replacing every vertex $v$ in $G$  by a copy $\dot{B}_v$  of  $\dot{B}$.
Let $K_2^3$ be the graph with two vertices and three parallel edges. We assume that $K_2^3$ is 3-connected by definition.
\\

From {\bf \ref{B}} we have:

\bs
\label{vRplB} 
Let $G$ be a cubic graph with possible parallel edges and 
$G' = G[\dot{B}]$. 
Then $v(G') = 9v(G)$, $\gamma (G') = 3v(G)$,    
$\kappa(G') = \kappa (G)$, and $G'$ is not cyclically 4-connected.
\es

The minimum cubic 3-connected graph provided by the above construction is $K_2^3[B]$.  
Obviously 
$v(K_2^3[B]) = 18$,  $\gamma  (K_2^3[B]) = 6$, and
$K_2^3[B]$ is obtained from two disjoint copies 
$(B', \{b'_1,b'_2,b'_3\})$ and $(B'', \{b''_1,b''_2,b''_3\})$ 
of $(B, \{b_1,b_2,b_3\})$ by adding three new edges 
$b'_ib''_i$, $i \in \{1,2,3\}$. 

Let $P^2_7$ be the Petersen (7,2)-graph.
Obviously $P^2_7$ is a cubic cyclically 4-connected graph with 14 vertices.
It can be checked that   
$\gamma (P^2_7) = 5 = \lceil v(P^2_7)/3\rceil  $ and
 $P^2_7$ is Hamiltonian.

Below (see {\bf \ref{Mk}}) we give constructions  that  for every $r \in \{0,1,2\}$ provide infinitely many cubic 3-connected and
cyclically 4-connected graphs $G_r$ such that 
$v(G_r) = r \bmod 3$,
$r \in \{0,2\} \Rightarrow \gamma (G_r) = \lceil v(G_r)/3 \rceil $,
and 
$r = 1 \Rightarrow \gamma (G_r) = \lfloor v(G_r)/3 \rfloor $.

Let $S$ be square $(t_1s_1t_2s_2t_1)$, $P$ be 4-vertex path $P = (q_1p_1p_2q_2)$.
Let $W$ be the graph obtained from disjoint $S$ and $P$ by identifying  $q_1$ with $s_1$ and $q_2$ with $s_2$.
Obviously $T = \{t_1,t_2,p_1,p_2\}$ is the set of degree two vertices in $W$.
\\

It is easy to prove the following.
\bs
\label{N,T}
Let $\dot{W} = (W, T)$ and $V \subseteq T$.
Then $\gamma (W - V) = 1$ if $V = \{p_1,p_2, t_i\}$ 
for some $i \in \{1,2\}$, and $\gamma (W - V) = 2$, otherwise.
\es

Let $k \ge 1$  be an integer,  $X = (x_0 \cdots x_{3k}) $ and 
$Y = (y_0 \cdots y_{3k})$ be two disjoint cycles, and 
$M^2_k = X \cup Y \cup \{x_0y_0, x_1y_1\} \cup 
\{x_iy_i: 1 \le i \le 3k-2, i = 1 \bmod 3\} \cup
\{x_iy_{i+1},x_{i+1}y_i: 2 \le i \le 3k-1,  i = 1 \bmod 3\}$.
%
Let 
$M^0_k = (M^2_k - \{x_0,y_0\}) \cup \{x_1x_{3k},y_1y_{3k}\}$, and $M^1_k = (M^2_k - \{x_0,y_0,x_1,y_1\}) \cup \{x_2x_{3k},y_2y_{3k}\}$.
Obviously 
$v(M^i_k)= i \bmod 3$.

\bs
\label{Mk} 
Each $M^i_k$ is a cubic cyclically 4-connected Hamiltonian graph and 
\\[0.5ex]
$(a0)$ $v(M^0_k)= 6k$ and $\gamma  (M^0_k)= 2k$,
\\[0.5ex]
$(a1)$ $v(M^1_k)= 6k - 2$ and $\gamma  (M^1_k)= 2k-1$,
and
\\[0.5ex]
$(a2)$ $v(M^2_k)= 6k +2$ and $\gamma  (M^0_k)= 2k+1$.
\es 

\bp (uses {\bf \ref{N,T}}). It is easy to see that each $M^i_k$,
$i \in \{0,1,2\}$,  is cyclically 4-connected and has 
a Hamiltonian cycle.
We prove $(a2)$. Claims $(a0)$ and $(a1)$ can be proved similarly. Obviously $v(M^2_k)= 6k +2$.

Since  $M^2_k$ is Hamiltonian, it has a dominating set with $2k +1$ vertices, and so $\gamma (M^2_k) \le 2k+1$.
Thus it is sufficient to show that if $D$ is a dominating set 
in $M^2_k$, then $|D| \ge 2k+1$.
We prove our claim by induction on $k$.
It is easy to check that our claim is true for $k \in  \{1, 2\}$.
So let $k \ge 3$.

Let $R_{3i+r}$ be the subgraph of $M^2_k$ induced by
the vertex subset
\\
$\{x_{3i+r}, x_{3i+r +1}, x_{3i+r+2}, y_{3i+r}, y_{3i+r+1},y_{3i+r+2} \}$, 
where  $i \in \{0, \ldots k-1\}$ and $r \in \{1,2\}$.
Then each  $R_{3i+r}$ is isomorphic to $W$ in 
{\bf \ref{N,T}} with 
$\{x_{3i+r+1}, y_{3i+r+1}\}$ corresponding to $\{s_1,s_2\}$,  
$V(R_{3i+r}) \cap  V(R_{3j+r}) = \emptyset $ for $i \ne j$,
$V(R_{3i+r}) \cap  \{x_{r-1},y_{r-1}\} = \emptyset $,
and 
$V(M^2_k) = \{x_{r-1},y_{r-1}\} \cup \{V(R_{3i+r}): 
i \in \{0, \ldots , k-1\}\}$.
Let $M = M^2_k$ and $R = R_{3i+1}$.
\\[1ex]
${\bf (p1)}$
Suppose that $M$ has a  minimum dominating set containing 
$Z = \{x_{3i+r},y_{3i+r}\}$ for some 
$i \in \{0, \ldots , k-1\}$ and $r \in \{2,3\}$. 
By symmetry of $M$, we can assume that $r = 2$.
Obviously $Z$ is a dominating set of $R$
and every degree two vertex in $R$ is adjacent to exactly one vertex in $M - R$.
Therefore $\gamma (M) = \gamma (M - R) + |Z|$.
Let $M' = (M - R) \cup \{x_{3i}x_{3i +4}, y_{3i}y_{3i +4}\}$.
Then $\gamma (M - R) \ge \gamma (M')$.
By the induction hypothesis, $\gamma (M') = 2k-1$.
Thus $\gamma (M) = \gamma (M - R) + |Z| \ge 
\gamma (M') + |Z| = (2k-1)+2 = 2k+1$.
\\[1ex]
${\bf (p2)}$ 
Suppose that $M$ has a  minimum dominating set  
$D$ containing one of the sets 
$\{x_{3i+r}, y_{3i+r+2}\}$, $\{y_{3i+r}, x_{3i+r+2}\}$,
$\{y_{3i+r}, y_{3i+r+2}\}$, $\{y_{3i+r}, y_{3i+r+2}\}$ for some
$i \in \{0, \ldots , k-1\}$  and $r \in \{1,2\}$.
By symmetry of $M$, we can assume that $D$ contains 
$\{x_{3i+1}, y_{3i+3}\}$ from $V(R)$.
If there is $z \in D \cap \{x_{3i+2}, y_{3i+2}\}$, then 
$D - z + x_{3i+3}$ is also 
a minimum dominating set of $M$.
Therefore we are done by {\bf (p1)}.
If $y_{3i+1} \in D$, then $D - y_{3i+1} + x_{3i}$ is also 
a minimum dominating set of $M$.
Thus we can assume that 
$D \cap V((R) = \{x_{3i+1}, y_{3i+3}\}$.
Then $D' = D \setminus \{x_{3i+1}, y_{3i+3}\}$ dominates 
$V(M) \setminus 
(\{x_{3i}, y_{3i+4}\} \cup V(R - x_{3i+3}))$.
Since  $D'$ dominates $x_{3i+3}$, clearly 
$x_{3i+4} \in D'$.
Let $M' = (M - R) \cup \{x_{3i}x_{3i+4}, y_{3i}y_{3i+4}\}$.
Then $M'$ is isomorphic to $M^2_{k-1}$ and since 
$x_{3i+4}$ dominates $\{x_{3i}, y_{3i+4}\} $, clearly $D'$ dominates $M'$. Therefore $|D'| \ge \gamma (M')$.
By the induction hypothesis, $\gamma (M') = 2k-1$.
Therefore 
$2k+1 \ge |D| = |D'|  + |\{x_{3i+1}, y_{3i+3}\}| =
(2k-1) + 2 = 2k+1$.
\\[1ex]
${\bf (p3)}$ 
Suppose that $M$ has a minimum dominating set $D$ containing one of the sets 
$\{x_{3i+r}, y_{3i+r+1}\}$, $\{y_{3i+r}, x_{3i+r+1}\}$
for some
$i \in \{0, \ldots , k-1\}$  and $r \in \{0,1\}$.
By symmetry of $M$, we can assume that $D$ contains 
$\{x_{3i+1}, y_{3i+2}\}$ from $V(R)$.
By ${\bf (p1)}$ and ${\bf (p2)}$, we can assume that
$D \cap \{x_{3i+2}, x_{3i+3}, x_{3i+4}, y_{3i+3}, y_{3i+4}\} 
= \emptyset $. 
Therefore $\{x_{3i+5}, y_{3i+5}\} \subseteq D$.
If $x_{3i+5}y_{3i+5} \not \in E(M)$, then we are done by 
${\bf (p1)}$. Therefore $x_{3i+5}y_{3i+5} \in E(M)$.
If $y_{3i+1} \in D$, then 
$D - y_{3i+1} + y_{3i}$ is also a minimum dominating set of $M$. Thus we can assume that 
$D \cap V(R) = \{x_{3i+1}, y_{3i+2}\}$.
Then $D' = D \setminus \{x_{3i+1}, y_{3i+2}\}$ dominates 
$V(M - x_{3i}) \setminus V(R)$.
Let $M'$ be as in ${\bf (p2)}$.
If $D' \cap \{x_{3i-1}, x_{3i}, y_{3i-1}\} \ne \emptyset $,
then $D'$ dominates $M'$, and we are done by the arguments similar to those in ${\bf (p2)}$.
If $D' \cap \{x_{3i-1}, x_{3i}, y_{3i-1}\} = \emptyset $, then
$y_{3i} \in D'$. By ${\bf (p2)}$, we can assume that
$D' \cap  \{x_{3i-2},  y_{3i-2}\} = \emptyset $.
Then $\{x_{3i-3},  y_{3i-3}\} \subseteq D'$.
Since $k \ge 3$, clearly $x_{3i-3} y_{3i-3} \not \in E(M)$.
Therefore we are done by ${\bf (p1)}$.
\\[1ex]
${\bf (p4)}$
Suppose that $M$ has a minimum dominating set $D$ that has exactly one vertex $z$ in $R_{3i+r}$ for some $i \in \{0, \ldots , k-1\}$ and $r \in \{1,2\}$. By symmetry of $M$, we can assume that
$r = 1$.
Then by {\bf \ref{N,T}}, $z \in \{x_{3i+3}, y_{3i+3}\}$. 
By symmetry of $M$, we can assume that 
$z = x_{3i+3}$, and so  by {\bf \ref{N,T}}, $y_{3i+4} \in D$.
Since $\{x_{3i+3}, y_{3i+4}\} \subseteq D$, we are done 
by ${\bf (p3)}$.
\\[1ex]
${\bf (p5)}$ Now suppose that for some 
$s \in \{0, \ldots , k-1\}$,
\\[0.5ex]
$(d1)$ a minimum dominating set $D$
contains exactly one of the  four sets 
$\{x_{3s+2}, x_{3s+3}\}$, $\{x_{3s+2}, y_{3s+3}\}$,
$\{x_{3s+3}, y_{3s+2}\}$,  and $\{y_{3s+2}, y_{3s+3}\}$. 

We can also assume by  ${\bf (p1)}$ and ${\bf (p2)}$ that 
\\[0.5ex]
$(d2)$ $D \cap \{x_{3s+1}, y_{3s+1}, x_{3s+4}, y_{3s+4}\}  
= \emptyset $. 

Then $(d1)$ and $(d2)$  hold for every $i \in \{0, \ldots , k-1\}$.
Hence $D \cap \{x_0, x_1, y_0, y_1\}  \ne \emptyset $ 
because  $D$ is a dominating set of $G$. Thus $|D| \ge 2k+1$.
\ep
\\

Let $N^r_k(i) = 
(M^2_k - \{x_{3i+1}x_{3i+2}, y_{3i+1}y_{3i}) \cup 
\{x_{3i+1}y_{3i}, y_{3i+1}x_{3i+2}\}$, where $1 < i < k$ and
$r \in \{0,1,2\}$.
One can also prove the following.

\bs
\label{Nk} 
Each $N^r_k(i)$ is a cubic 3-connected $($but not cyclically 
4-connected$)$ Hamiltonian graph and 
\\[0.5ex]
$(a0)$ $v(N^0_k(i))= 6k$ and $\gamma  (N^0_k(i))= 2k$,
\\[0.5ex]
$(a1)$ $v(N^1_k(i))= 6k - 2$ and $\gamma  (N^1_k(i))= 2k-1$,
and
\\[0.5ex]
$(a2)$ $v(N^2_k(i))= 6k +2$ and $\gamma  (N^0_k(i))= 2k+1$.
\es

We believe that the following is true.
\bs {\em \bf Conjecture}
\label{conjecture}
~ Let $G$ be a cubic 3-connected graph.
If $v(G) \ne 1 \bmod 3$, then $\gamma (G) \le \lceil v(G)/3 \rceil $.
If $v(G) = 1 \bmod 3$, then $\gamma (G) \le \lfloor v(G)/3 \rfloor$.
\es 

From {\bf \ref{vRplB}}, {\bf \ref{Mk}}, and {\bf \ref{Nk}} 
it follows that Conjecture {\bf \ref{conjecture}} is best possible
for both 3-connected and cyclically 4-connected cubic graphs. 

From the results in \cite{Kpack} it follows that  
if $G$ is a Hamiltonian cubic graph with $v(G) = 1 \bmod 3$, then $\gamma (G) \le \lfloor v(G)/3 \rfloor$.
Therefore Conjecture 
{\bf \ref{conjecture}} is true for Hamiltonian cubic graphs.





\begin{thebibliography}{99}

\bibitem{Ar} V.I. Arnautov, 
Estimation of the external stability number  of a graph by means of the minimum degree of vertices,
{\em Prikl. Math. i Program.} 10 (1974) 3--8.

\bibitem{BC} B. Bollob\'as and E.J. Cockayne,
Graph-theoretical parameters concerning domination, independence, and irredundance,
{\em J. Graph Theory} 3 (1979) 241--249.

\bibitem{B} M. Blank,
An estimate of the external stability of a graph without pendant
vertices,
{\em Prikl. Math. i Program.} 10 (1973) 3--11.


\bibitem{CH}E.J. Cockayne and S.T. Hedetniemi,
Towards a theory of domination in graphs,
{\em Network} 7 (1977) 247--261.

\bibitem{D} R. Diestel, {\em Graph Theory}, Springer, 2005.

\bibitem{HHS} T.W. Haynes, S.T. Hedetniemi, and P.J. Slater, {\em Fundamentals of domination in graphs}, Marcel Dekker, Inc., 1998.

\bibitem{KPS} K. Kawarabayashi, M. Plummer, and A. Saito,
Domination in a graphs with a 2-factor,
{\em J. Graph Theory} 52 (2006) 1--6.

\bibitem{Kpack} A. Kelmans, On packings in cubic graphs,
submitted.


\bibitem{KS} A.V. Kostochka and B.V. Stodolsky,
On domination in  graphs with minimum degree two,
{\em Discrete Mathematics} 304 (2005) 749--762.

\bibitem{MS} W. McCuaig and B. Shepherd,
Domination in connected cubic graphs,
{\em J. Graph Theory} 13 (1989) 45--50.

\bibitem{O} O.  Ore,  {\em Theory of Graphs}, Amer. Math. Soc. Colloq. Publ. 38, 1962.

\bibitem{R} B. Reed, Paths, stars, and the number three,
{\em Combin. Probab. Comput.} 5 (1996) 277-295.



\end{thebibliography}
\end{document}